\documentclass[runningheads]{llncs}

\usepackage{epsfig}

\usepackage{amsmath, amssymb}

\usepackage{algorithm}
\usepackage{algorithmic}

\usepackage{enumerate}

\begin{document}

\title{Nonparametric Markovian Learning of Triggering Kernels for Mutually Exciting and Mutually Inhibiting Multivariate Hawkes Processes}

\titlerunning{Markovian Estimation of Mutually Interacting Processes}

\author{Remi Lemonnier\inst{1,2} \and Nicolas Vayatis\inst{1}}

\institute{CMLA-ENS Cachan,Cachan,France  \and 1000Mercis,Paris,France}

\maketitle

\begin{abstract} 
In this paper, we address the problem of fitting multivariate Hawkes processes to potentially large-scale data in a setting where series of events are not only mutually-exciting but can also exhibit inhibitive patterns. We focus on nonparametric learning and propose a novel algorithm called MEMIP (Markovian Estimation of Mutually Interacting Processes) that makes use of polynomial approximation theory and self-concordant analysis in order to learn both triggering kernels and base intensities of events. Moreover, considering that N historical observations are available, the algorithm performs log-likelihood maximization in $O(N)$ operations, while the complexity of non-Markovian methods is in $O(N^{2})$. Numerical experiments on simulated data, as well as real-world data, show that our method enjoys improved prediction  performance when compared to state-of-the art methods like MMEL and exponential kernels.
\end{abstract} 

\section{Introduction}

Multivariate Hawkes processes are a class of multivariate point processes which are often used to model counting processes where physicals events rate of occurrence usually depend on past occurences of many other events. This is typically the case for earthquakes aftershocks (\cite{ogata1988statistical}) and financial trade orders on marketplace (\cite{errais2007pricing}, \cite{bauwens2009modelling},\cite{bacry2013modelling}, \cite{alfonsi2014dynamic}), but also in other fields such as crime prediction (\cite{mohler2011self}), genome analysis (\cite{reynaud2010adaptive}) and more recently for modeling social interactions (\cite{blundell2012modelling}, \cite{zhou2013learning}). Multivariate Hawkes processes are fairly well-known from a probabilistic point of view : their Poisson cluster representation was outlined by the seminal paper of Hawkes and Oakes (\cite{hawkes1974cluster}), stability conditions and sample path large deviations principles were derived in a sequence of papers by Bremaud and Massoulie (see e.g \cite{bremaud1996stability}). In the unidimensional case,  Ogata \cite{ogata1978asymptotic} showed that the log-likelihood estimator enjoys usual convergence properties under mild regularity conditions. However, in practical applications, estimation of the triggering kernels $g_{uv}$ has always been a difficult task. First, because Hawkes log-likelihood contains the logarithm of the weighted sum of triggering kernels, most of the aforementioned papers made the choice of fixing triggering kernels up to a normalization factor in order to ensure concavity, that is $g_{uv}=c_{uv} \cdot g$. Secondly, when computational efficiency is an issue, the dependency of the stochastic rate at a given time on all the past occurences implies quadratic complexity in the number of occurences for tasks like log-likelihood computation. This issue has often been tackled by choosing memoryless exponential triggering kernels, but the actual dynamics of kernels strongly depends on the field of application: price impacts of a given trade (\cite{bacry2012non}) and process of views of Youtube videos (\cite{crane2008robust}) were shown to be better described by slowly decaying power-law kernels whereas for DNA sequence modelization (\cite{reynaud2010adaptive}) kernels are known to have bounded support. Thus, it is highly desirable to estimate triggering kernels in a data-driven way instead of assuming a given parametric form. Nonparametric estimation has been successfully addressed  for unidimensional  (\cite{reynaud2010adaptive},\cite {lewis2011nonparametric}), and symmetric bidimensional (\cite{bacry2012non}) Hawkes processes. In the case where triggering kernels are known to sparsely decompose over a dictionary of basis functions of bounded support (e.g for neuron spikes interactions), a LASSO-based algorithm with provable guarantees was derived in \cite{hansen2012lasso}. 

Recently, combining majorization-minimization techniques with resolution of a Euler-Lagrange equation,  Zhou, Zha and Song (\cite{zhou2013learning}) proposed what is to our knowledge the first nonparametric learning algorithm for general multivariate Hawkes processes. But although this work constitutes a significant improvement over existing parametric methods, it still relies on several assumptions. First, interactions between events are assumed to be "mutually-exciting", i.e $g_{uu'}$ are non-negative for all $u,u'$. We nevertheless argue that in real-world settings, there is no reason to think that interactions beween events are only mutually-exciting. Secondly, the background rates $\mu_{u}$ are assumed to be constant. While this is a common assumption for multivariate Hawkes processes, it was shown by \cite{lewis2011self} that estimating $\mu_{u}(t)$ from the data could lead to significant improvement. To address these different issues, we construct a novel algorithm MEMIP (Markovian Estimation of Mutually Interacting Processes) based on polynomial approximation of a mapping of the triggering kernels to $[0,1]$. Our method does not assume non-negativity on triggering kernels and is able to estimate time-dependent background rate on a data-driven way. Moreover, by constructing a markovian and linear estimator, we carry the more appealing properties of the most widely used parametric setting, where triggering kernels are fixed to exponentials up to a normalization factor : concavity of the log-likelihood that ensures global convergence of the estimator, and $O(N)$ log-likelihood calculation in a single pass through the data. While giving a concave formulation of the exact log-likelihood that can be maximized by multiple optimization techniques, we propose an algorithm based on maximisation of a self-concordant approximation that is shown to outperform state-of-the-art methods on both simulated and real-world data sets.

The paper is organized as follows. In Section 2, we formally define multivariate Hawkes processes as well as the associated log-likelihood maximization problem. In section 3, we decompose the log-likelihood on a basis of memoryless triggering kernels. Through Section 4, we develop two novel algorithms for exact as well as fast approximate maximization of the log-likelihood, analyze their complexity and show numerical convergence results based on the properties of self-concordant functions. In section 5, we show that MEMIP significantly improves over state of the art on both synthetic and real world data sets for the tasks of predicting future events as well as estimating underlying dynamics of the Hawkes process. 

\section{Setup and Notations}

\subsection{Model Description and Notation}

We consider a multivariate Hawkes process, that is a $d$-dimensional counting process $N(t)= \{N^u(t)~:~u=1, \ldots, d\}$ for which the  rate of occurence of each component $N^u(t)$ is defined by:
\begin{equation}
\lambda_{u}(t)=\left( \mu_{u}(t)+\sum_{v \in [1...d]} \sum_{t_{v}<t} g_{uv}(t-t_{v}) \right)_{+}~, \quad \forall u =1, \ldots, d
\label{eqn:linear Hawkes}
\end{equation}
\noindent where $\mu_{u}(t)$ is the natural rate of occurence of events along  dimension $u$. Note that the occurence of a given event affects stochastic rates of occurence of every dimension. With an empty history, events of type $u$ will occur as if they were drawn from a non-homogeneous Poisson process of rate $\mu_{u}(t)$. The kernel function evaluation $g_{uv}(t-t_{v})$ quantifies the change in the rate of occurence of event $u$ at time $t$ caused by the realization of event $v$ at time $t_{v}$. Following the intuition, we can characterize three situations depending on the values taken by the kernel function at a given time lapse $s$:
\begin{itemize} 
\item \emph{Excitation} corresponds to the case where we have $g_{uv}(s) >0$,  {\em i.e.} an event of type $v$ is more likely to occur if an event of type $u$ has occured at a time distance of $s$. 
\item \emph{Independence} is observed when  $g_{uv}(s) = 0$,  meaning  that the realization of an event of type $u$ has no effect on the rate of occurence of an event of type $v$ at time distance $s$.
\item  \emph{Inhibition} takes place when $g_{uv}(s)<0$,  {\em i.e.} an event of type $v$ is less likely to occur if an event of type $u$ occured at time distance $s$.
\end{itemize}
Such processes can be seen as a generalization over the common definition of multivariate Hawkes process where the kernels $g_{uv}$ are non-negative and the componentwise background rate $\mu_{u}$ is often taken constant. 

\subsection{Log-Likelihood of Multivariate Hawkes Processes}

\subsubsection{Input Observations.} We define a \emph{realization} $h$ of a multivariate point process by the triplet $T_{h}^{-},T_{h}^{+},(t_{i}^{h},u_{i}^{h})_{i \in [1...n_{h}]}$, where $T_{h}^{-}$ and $T_{h}^{+}$ are respectively the beginning and the end of the observation period, and  $(t_{i}^{h},u_{i}^{h})$, for $i \in [1...n_{h}]$, is the sequence of the $n_{h}$ events occuring during this period. In the rest of the paper, we will assume we are given $n$ i.i.d realizations of a multivariate Hawkes process. Without loss of generality, we will assume $\min_{h}(T_{h}^{-})=0$ and take $T=\max_{h}(T_{h}^{+})$.

\subsubsection{Expression of the Log-Likelihood.} We first set $\Lambda = \{\lambda_u~:~u =1, \ldots, d\}$. For a general multivariate point process, the log-likelihood of the whole dataset $\mathcal{H}$ is given by ({\em e.g.} \cite{daley2007introduction}):
\begin{equation}
\mathcal{L}(\Lambda, \mathcal{H}) = \sum_{u=1}^d \sum_{h\in\mathcal{H}}\int_{T_{h}^{-}}^{T_{h}^{+}}\ln(\lambda_{u}(s))dN^{u}_h(s)-\sum_{u=1}^d \sum_{h\in\mathcal{H}}\int_{T_{h}^{-}}^{T_{h}^{+}} \lambda_u(s) ds 
\end{equation}
\noindent where $\int f(s))dN^{u}_h(s)=\sum_{i=1}^{n_{h}} f(t_{i}^{h})1\left\{u_{i}^{h}=u\right\}$. In the case of a linear Hawkes process (\ref{eqn:linear Hawkes}), we introduce $\Lambda = (M, G)$ where  $M=\{\mu_u~:~u =1, \ldots, d\}$ and $G = \{ g_{u,v}~:~u, v =1, \ldots, d \}$ and the log-likelihood can be rewritten as:
\begin{eqnarray}
\mathcal{L}(M, G, \mathcal{H}) = \sum_{h\in\mathcal{H}} \sum_{i=1}^{n_{h}} \ln\bigg(\mu_{u_{i}^{h}}(t_{i}^{h})+\sum_{j~:~t_{j}^{h}<t_{i}^{h}}  g_{u_{j}^{h},u_{i}^{h}}(t_{i}^{h}-t_{j}^{h})\bigg)\nonumber\\ 
-\sum_{u=1}^{d} \sum_{h\in\mathcal{H}} \int_{T_{h}^{-}}^{T_{h}^{+}} \bigg(\mu_{u}(s)+\sum_{j=1}^{n_{h}}1\left\{u_{j}^{h}=u\right\}\ g_{u,u_{j}}(s-t_{j})\bigg)_{+} ds
\label{eqn:Log-likelihood-pos}
\end{eqnarray}
Depending on the parametrization of triggering kernels $g_{uv}$, this log-likelihood may or may not be concave. For instance, in the widely used setting where the background rates $\mu_{u}$ are constant and the kernels $g_{uv}$ are non-negative and fixed up to the normalization factor $\nu_{uv}$, the log-likelihood is concave and can be relatively easily maximized. However, even for the simple case of nonnegative exponential kernels $g_{uv}(t)=\nu_{uv}\exp(-\alpha_{j}t)$ where $\nu_{uv}\ge 0$ the product term $\nu_{uv}\exp(-\alpha_{v}t)$ makes the log-likelihood not concave with respect to $\alpha_{v}$. Therefore, global convergence of maximization methods is not guaranteed anymore.

\section{Approximations of Multivariate Hawkes Processes on a Basis of Exponential Triggering Kernels}

\subsection{A $K$-approximation of the Multivariate Hawkes Process} For a given multivariate Hawkes process $\Lambda = (M, G)$, we consider finite approximations of the components of the  rates of occurence $\mu_u$ and $g_{uv}$. We first introduce the following functions:
\[
\forall y\in [-\ln(T)/\alpha, 1], \quad \nu_{u} (y) = \mu_{u}(-\ln(y) / \alpha)  \quad  \text{and} \quad  f_{uv} (y) = g_{uv}(-\ln(y) / \alpha)
\] 
and we  use Bernstein-type polynomial approximations of order $K$ for $\nu_u$ and  $f_{uv}$: there exist coefficients $X^K_{uv, k}$ such that
\[
\forall y\in [-\ln(T)/\alpha, 1], \quad  \widehat{\nu}^K(y) = \sum_{k=0}^K X^K_{u0, k}y^k \quad \text{and} \quad  \widehat{f}^K_{uv} (y) =  \sum_{k=0}^K X^K_{uv, k}y^k ~.
\]
These polynomial approximations are known to converge with a polynomial rate for smooth functions (with first $r$ derivatives continuously differentiable) and geometric rate for analytic functions (see below). The $K$-aproximation considered in this paper relies on a simple change of variable in the Bernstein approximations by setting: $y=\exp(-\alpha t)$. 
We can now introduce the linear approximation of a multivariate Hawkes process with exponential kernels:
\[
\forall t\in [0,T], \quad  \widehat{\mu}^K(t) = \sum_{k=0}^K X^K_{u0, k}\exp(-k\alpha t) \quad \text{and} \quad  \widehat{g}^K_{uv} (t) =  \sum_{k=0}^K X^K_{uv, k}\exp(-k\alpha t)~.
\]

Classical arguments from approximation theory (\cite{bernstein1912ordre} and \cite{cheney1966introduction}) lead to the following proposition.
\begin{proposition} \label{prop:approx-rates}
For any function $\Psi$ defined over $[0,T]$, we consider the supremum norm $\left |\left| \Psi \right|\right|_{T,\infty} = \sup_{t \in [0,T]} |\Psi(t)|$. The $K$-approximations  ($\widehat{\mu}^{K}_{u})_{K\ge 1}$ and ($\widehat{g}^{K}_{uv})_{K\ge 1}$ converge in supremum norm towards true functions $\mu_{u}$ and $g_{uv}$ at the following rates:
\begin{enumerate}
\item if $\mu_{u}$ is $C^{r}$, $\left |\left|\mu_{u}(t)-\widehat{\mu}^{K}_{u}(t)\right|\right|_{\infty}^{T}=O(1/K^{r})$ 
\item if $\mu_{u}$ is analytic, $\left |\left|\mu_{u}(t)-\widehat{\mu}^{K}_{u}(t)\right|\right|_{\infty}^{T}=O(\exp(-K))$ 
\item if $g_{uv}$ is $C^{r}$, $\left |\left|g_{uv}(t)-\widehat{g}^{K}_{uv}(t)\right|\right|_{\infty}^{T}=O(1/K^{r})$ 
\item if $g_{uv}$ is analytic, $\left |\left|g_{uv}(t)-\widehat{g}^{K}_{uv}(t)\right|\right|_{\infty}^{T}=O(\exp(-K))$.
\end{enumerate}
\end{proposition}

Another property of the approximated multivariate Hawkes process is the Markov property of the counting process. We set  $\widehat{N}^{K}(t)$ the  $d$-dimensional Hawkes process uniquely defined by $\widehat{\lambda}^{K} = (\widehat{\mu}^{K}_{u}, \widehat{g}^{K}_{uv})_{u, v}$.
\begin{proposition}
Assume that the empirical estimate $\widehat{N}^{K}(t)$ of the multivariate Hawkes process is obtained after $i.i.d.$ realizations of $N(t)$ over the time interval $[0,T]$. There exists $(\widehat{\ell}^{0},\widehat{\ell}^{1}, \ldots, \widehat{\ell}^{K})$ such that:
 \[
\forall u\in \{1, \ldots, d\}~, \quad\widehat{\lambda}^{K}(t)=\sum_{k=0}^{K} \left(\widehat{\ell}^{k}(t) \right)_{+}
\]
and ($\widehat{N}^{K}(t),\widehat{\ell}^{0}(t),\widehat{\ell}^{1}(t), \ldots, \widehat{\ell}^{K}(t))$ is a Markov Process on $\mathbb{N}^{d} \times \mathbb{R}^{d(K+1)}$. 
\end{proposition}

The proof results from the following decomposition of each occurrence rate in the approximation: $\forall u\ge 1$,
\begin{multline}
\widehat{\lambda}^{K}_{u}(t)= \Bigg(X^{K}_{u0,0}+\sum_{k=1}^{K} \bigg(X^{K}_{u0,k} \exp(-k\alpha t) +  \sum_{v~:~t_{v}<t} X^{K}_{uv,(k-1)} \exp(-k\alpha (t-t_{v}))\bigg)\nonumber\\+ \sum_{v~:~t_{v}<t} X^{K}_{uv, K} \exp(-(K+1)\alpha (t-t_{v}))\Bigg)_{+}
\end{multline}

Markov property is then a direct consequence of the dynamics of the functions $\widehat{\ell}^{k}_{u}(t)$ : they decay at rate $\exp (-k \alpha t)$ and jump by $X^{K}_{uv,(k-1)}$ whenever an event of type $v$ occurs. As they entirely determine the stochastic rate which determines the conditional probability distribution of $\widehat{N}^{K}(t)$, the conditional probability distribution of future states of the process $(\widehat{N}^{K}(t),\widehat{\ell}^{0}(t),\widehat{\ell}^{1}(t),...\widehat{\ell}^{K}(t))$ is uniquely determined by the present state.

\subsection{A New Decomposition of the Log-Likelihood}

The algorithms proposed in this paper rely on a novel expression of the log-likelihood over a basis of triggering kernels. We use exponential excitation functions to account for nonlinearity but our algorithms benefit from the properties of linear approximations. Based on the expression of the log-likelihood for general linear multivariate Hawkes process (\ref{eqn:Log-likelihood-pos}), we introduce the following notation to discover the specific expression for the $K$-approximation based on exponential triggering functions: $\forall u, v = 1, \ldots, d$, $\forall k=1, \ldots, K$, $\forall h\in \mathcal{H}$,  $\forall i=1, \ldots, n_h$,
\begin{align}
A^{K, h,i}_{uv,k} &= \sum_{j~:~t_{j}^{h}<t_{i}^{h}}1\left\{u_{i}^{h}=v, u_{j}^{h}=u\right\}   \exp \bigl(-(k+1\left\{u>0\right\})\alpha (t_{i}^{h}-t_{j}^{h})\bigr)\\
B^{K,h}_{0v,k}(s)&=\exp(-k\alpha s) \\
B^{K,h}_{uv,k}(s)& =\sum_{j ~:~t_{j}^{h}<s} 1\left\{u_{j}^{h}=v\right\} \exp(-(k+1)\alpha (s-t_{j}^{h}))
\end{align}
The key expression of the approximate log-likelihood can then be derived by plugging-in the previous notations and replacing the intrinsic parameters $(M, G)$ by the linear coefficients $X^K$: 
\begin {eqnarray}
\mathcal{L}^{K}(X^{K}, \mathcal{H})=
\sum_{h\in \mathcal{H}} \sum_{i=1}^{n_{h}} \ln(A^{K,h,i}X^{K})
- \sum_{h\in \mathcal{H}}\int_{0}^{T_{h}} \bigg (\sum_{i=1}^{n_{h}} B^{K,h}(s)X^{K}\bigg)_{+}ds 
\label{eqn:Log-likelihood_P_vectorized}
\end{eqnarray}

Note that the dependance of $\mathcal{L}^{K}$ on the history $\mathcal{H}$ is entirely expressed by vectors $(A^{K, h,i})_{h \in \mathcal{H},i \in [1...n_{h}]}$ and $(B^{K,h}(s))_{h \in \mathcal{H},s\in[0,T]}$. An important feature of the approximate log-likelihood expressed in the parameter space defined by linear decompositions onto bases of exponential triggering kernels is given in the following proposition.

\begin{proposition} 
The function $X \rightarrow \mathcal{L}^{K}(X, \mathcal{H})$ is concave.
\end{proposition}

From there, we have a complete roadmap for the design of algorithms estimating the parameters of multidimensional Hawkes processes: the last propostion indicates that a proxy of the log-likelihood (\ref{eqn:Log-likelihood-pos}) can  be globally maximized with tools of convex analysis. Moreover, thanks to the approximation rates of convergence (Proposition \ref{prop:approx-rates}), triggering kernels can be accurately estimated for large $K$ 
through maximization of the new objective (\ref{eqn:Log-likelihood_P_vectorized}). Finally, the Markov property is an important feature that will allow us to construct the vectors $(A^{K,h,i})$ and $(B^{K,h})$ with linear complexity.

\section{Markovian Algorithms for the Estimation of Triggering Kernels}

Computational tractability of algorithms on large data sets depends on the algorithmic complexity in the dominating dimensions of the problem. For realizations of multivariate Hawkes processes, dominating dimensions are almost always the total number of events $N=\sum_{h \in \mathcal{H}} n_{h}$ and the time of observation $T$. Indeed, it would be unrealistic to try to learn $d^{2}$ nonparametric functions in an infinite dimensional space with only $N$ observations without the condition $N \gg d^{2}$. In the rest of the paper, we will therefore focus on constructing two algorithms with no more than linear complexity in $N$ and $T$.

\subsection{Exact Maximization of the Approximated Log-Likelihood }

Vectors $(A^{K,h,i})_{h \in \mathcal{H},i \in [1...n_{h}]}$ and $(B^{K,h}(s))_{h \in \mathcal{H},s\in[0,T]}$ can be constructed in a single pass through the data by \textbf{Algorithm 1}.

\begin{algorithm}
\caption{Algorithm for construction of vectors $(A^{K,h,i})$ and $(B^{K,h}(s))$}
\begin{algorithmic}
\STATE {Initialize $i=0$ and fix a time step $dt$}
\FORALL{$h$}
\STATE{Initialize $(C_{uv}^{k}=0)_{u \geq 1, v \geq 1}$ ; $t=T_{h}^{-}$ ; $(D_{uv}^{k}(T_{h}^{-})=1_{\left\{u=0\right\}})_{u \geq 0, v \geq 1}$}
\WHILE{$t < T_{h}^{+}$}
\STATE{$t \leftarrow t+\delta t=\min(t+dt,t_{i})$}
\FORALL{k,u,v}
\STATE{$C_{uv}^{k} \leftarrow C_{uv}^{k} \exp(-(k+1\left\{u>0\right\} \alpha \delta t)$, $D_{uv}^{k} \leftarrow D_{uv}^{k} \exp(-(k+1\left\{u>0\right\} \alpha \delta t)$}
\STATE{$B^{K,h}_{uv,k}(t) \leftarrow D_{uv}^{k}$}
\ENDFOR
\IF{$t=t_{i}$}
\FORALL{k,u}
\STATE{$A^{K,h,i}_{uv,k} \leftarrow C_{uu_{i}}^{k}$}
\ENDFOR
\FORALL{k,v}
\STATE{$C_{u_{i}v}^{k} \leftarrow C_{u_{i}v}^{k} + 1$, $D_{u_{i}v}^{k} \leftarrow D_{u_{i}v}^{k}+ 1$}
\ENDFOR
\STATE{$i \leftarrow i+1$}
\ENDIF
\ENDWHILE
\ENDFOR
\end{algorithmic}
\end{algorithm}

\medskip

\noindent{\bf Complexity of  \textbf{Algorithm 1}.} With $M=T/dt$ the number of discretizations steps, construction of vectors $(A^{K,h,i})$ and $(B^{K,h}(s))$ has thus a complexity of $O(N+M)$. As each log-likelihood evaluation (\ref {eqn:Log-likelihood_P_vectorized}) requires $2N+M$ scalar products computations, various optimization techniques can be used to find the global maximum of $X \rightarrow \mathcal{L}^K(X, \mathcal{H})$ in $O(N+M)$ operations. On the contrary, a nonmarkovian estimator, even linear, would need at each time $t$ to compute the values of triggering kernels between current time and all preceding occurence times, thus leading to a $O(\sum_{h}n_{h}^{2})$ complexity. This construction is thus very often the bottleneck of the whole maximization procedure.

\subsection{Relaxed Version of the Log-Likelihood}

While the previous paragraph exposes a fully tractable method to estimate the triggering kernels for potentially large data sets, we now develop an approximate algorithm called MEMIP, for Markovian Estimation of Mutually Interacting Processes, that leads to a substantial speed-up, as well as theoretical  guarantees  in terms of efficiency. For this purpose, we approximate the log-likelihood $\mathcal{L}^K(M, G, \mathcal{H})$ by dropping the positive part in log-likelihood (\ref{eqn:Log-likelihood-pos}), {\em i.e.}
\begin{eqnarray}
\widetilde{\mathcal{L}}^K(M, G, \mathcal{H})= \sum_{h\in\mathcal{H}} \Bigg( \sum_{i=1}^{n_{h}} \ln\bigg(\mu_{u_{i}^{h}}(t_{i}^{h})+\sum_{j~:~t_{j}^{h}<t_{i}^{h}}  g_{u_{j}^{h},u_{i}^{h}}(t_{i}^{h}-t_{j}^{h})\bigg) \nonumber\\  -
\sum_{u=1}^{d} \int_{T_{h}^{-}}^{T_{h}^{+}} \bigg(\mu_{u}(s)+\sum_{j=1}^{n_{h}}1\left\{u_{j}^{h}=u\right\} g_{u,u_{j}}(s-t_{j})\bigg) ds \Bigg)
\end{eqnarray}
which can be rewritten:
\begin{equation}
 \widehat{\mathcal{L}}^K(X^{K}, \mathcal{H} )=\sum_{h\in\mathcal{H}} \bigg(\sum_{i=1}^{n_{h}} \ln(A^{K,h,i}X^{K})\bigg)- \widehat {B}^{K}X^{K}
\label{eqn:LLhat}
\end{equation} 
where $\displaystyle \widehat{B}^{K}_{uv,k}= \sum_{h\in\mathcal{H}} \sum_{j=1}^{n_{h}} 1\left\{u_{j}^{h}=v\right\} \int_{T_{h}^{-}}^{T_{h}^{+}}  \exp(-k\alpha (s-t_{j}^{h}))$.

\medskip

Although $ \widehat{\mathcal{L}}^K(X, \mathcal{H} )$ is an upper bound of the actual log-likelihood and it is not clear at first sight why its maximization should lead to large values of $\mathcal{L}^K(X, \mathcal{H} )$, we point out that the difference $\widehat{\mathcal{L}}^K(X, \mathcal{H} )-\mathcal{L}^K(X, \mathcal{H} )$ is only caused by intervals where there exists $u \in[1...d]$ such that $\widehat{\lambda}^{K}_{u}(t) = 0 $. But maximizers of $\widehat{\mathcal{L}}^K(X, \mathcal{H} )$ are very unlikely to exhibit wide range of negative values in their triggering kernels because any single event realization with a predicted nonpositive stochastic rate yields $\widehat{\mathcal{L}}^K(X, \mathcal{H} ) = - \infty$. Therefore, we assume we can rely on this approximation in order to construct fast algorithms.

\subsection{MEMIP: a Learning Algorithm for Fast Log-Likelihood Estimation}

Since the gradient and the hessian matrix of $X \mapsto \widehat{\mathcal{L}}^K(X, \mathcal{H})$ can be computed analytically and their size does not depend on $N$, we derive the proposed algorithm MEMIP on the base of successive Newton optimizations. In the following, we denote by $\text{NewtonArgMax}(f,x_{0})$ the result of a Newton maximization of function $f$ with starting point $x_{0}$ using a classical backtracking linesearch method. The main idea is to construct recursively a sequence $(\widehat{X^{1}}...\widehat{X^{K}})$ of maximizers of functions $( \widehat{\mathcal{L}}^k)_{k \in [1...K]}$ by using $\text{NewtonArgMax}(\widehat{\mathcal{L}}^{k-1}, \widehat{W}^{k-1})$ as the starting point $\widehat{W}^{k}$ of maximization of $ \widehat{\mathcal{L}}^{k}$.
\begin{algorithm}
\caption{Algorithm (MEMIP) for learning background rates and triggering kernels of a multivariate Hawkes process}
\begin{algorithmic}
\INPUT{Mapping parameter $\alpha>0$, maximal polynomial degree $K$, starting point $\widehat{W}^{1} \in \mathbb{R}^{d(d+1)}$}
\STATE{Construct $(A^{K,h,i})$ and $B^{K}$ according to $O(N)$ modified version of {\bf Algorithm 1}}
\STATE{$\widehat{X}^{1} \leftarrow \text{NewtonArgMax}(\widehat{\mathcal{L}}^1, \widehat{W}^{1})$}
\FOR{$k \in [2...K]$}
\STATE{$\widehat{W}^{k} = 0$}
\FOR{$j \in [1...k-1], u \in[1...d], v\in[0...d]$}
\STATE{$\widehat{W}^{k}_{uv,j}=\widehat{X}_{uv,j}^{k-1}$}
\ENDFOR
\STATE{$\widehat{X}^{k} \leftarrow \text{NewtonArgMax}(\widehat{\mathcal{L}}^k, \widehat{W}^{k})$}
\ENDFOR
\end{algorithmic}
\end{algorithm}
From the estimated sequence $(\widehat{X}^{1}...\widehat{X}^{K})$, the best value of $k$ can be estimated by cross-validation or various other model selection techniques. Interestingly, $A^{k,h,i}=(A^{K,h,i}_{\bullet,j})_{j \in [1...k]}$ and $B^{k}=(B^{K}_{\bullet,j})_{j \in [1...k]}$ such that only $(A^{K,h,i})_{h \in \mathcal{H},i \in [1...n_{h}]}$ and $B^{K}$ need to be computed.

\medskip

\noindent{\bf Complexity of  \textbf{Algorithm 2}.} We obtain two substantial computational speed-ups compared to exact log-likelihood maximization. First, time discretization is no longer needed for the construction of $B^{K}$. Thus, vectors $(A^{K,h,i})$ and $B^{K}$ can be constructed with the same procedure than \textbf{Algorithm 1} except that updates are made only on time occurence of events. Therefore, construction complexity is $O(N)$. Similarily, approximate log-likelihood evaluations are also of complexity $O(N)$. Secondly, the approximate log-likelihood is separable by type of event $u$ : $\widehat{\mathcal{L}}^K = \sum_{u=1}^{d} \widehat{\mathcal{L}}^K_{u}$ where $\widehat{\mathcal{L}}^K_{u}$ only depends on background rate $\mu_{u}$ and triggering kernels $(g_{uv})_{v \in [1...d]}$. Maximization can thus be parallelized across the different dimensions. Note that because of the Hessian inversion at each Newton step, complexity in d of maximization of $\widehat{\mathcal{L}}^K_{u}$ is $O(d^{3})$ for any $u$, which yields a $O(d^{4})$ overall complexity. In cases where $N \gg d^{2}$ but $d^{4}>N$, the use of quasi-Newton methods might therefore be preferable.

\subsection{Self-Concordance Property and Numerical Convergence of MEMIP}

Problem (\ref{eqn:LLhat}) can be solved by various optimisation techniques. {\bf Algorithm 2} is actually based on the concept of \emph{self-concordance}  (\cite{nesterov1994interior}) that we apply to function $ X \mapsto-\widehat{\mathcal{L}}^k(X, \mathcal{H})$. Self-concordant functions are, along with strongly-convex functions with Lipschitz-continuous Hessian matrices, a very important class of functions for which nonasymptotic upper bounds of the number of Newton steps necessary to reach precision $\epsilon$ is known.  More specifically, the following property holds:

\begin{proposition} Starting from a d(d+1)-dimensional vector $\widehat{W}^{1}$, MEMIP constructs a sequence of $K$ estimates $(\widehat{X}^{1}...\widehat{X}^{K})$ verifying for any $k \in [1...K]$, $|\widehat{\mathcal{L}}^{k}(\widehat{X}^{k},\mathcal{H})-\sup_{X}(\widehat{\mathcal{L}}_{k}(X,\mathcal{H}))| \leq \epsilon$ in at most $C \big(\sup_{X}(\widehat{\mathcal{L}}_{K}(X,\mathcal{H}))-\widehat{\mathcal{L}}_{1}(\widehat{W}^{1},\mathcal{H})\big)+K (\log_{2}\log_{2}(1/\epsilon)+C\epsilon )$ Newton iterations.
\end{proposition}

\begin{lemma}
Using Newton method with backtracking line search from a starting point $x_{0} \in \mathbf{R}^{d}$, there exists $C>0$ depending only on the line search parameters such that the total number of Newton iterations needed to minimize a self-concordant function $f$ up to a precision $\epsilon$ is upper bounded by $C(\sup(f)-f(x_{0}))+\log_{2}\log_{2}(\frac{1}{\epsilon})$.
\end{lemma}

\noindent {\bf Proof of Proposition 4.} Self-concordance of functions ($-\widehat{\mathcal{L}}_{k})_{k\in[1...K]}$ is a direct consequence of self-concordance on $\mathbf{R}^{*}_{+}$ of $f:x \mapsto -\ln(x)$ and affine invariance properties of self-concordant functions. By applying the aforementioned lemma to function $-\widehat{\mathcal{L}}_{k}$ and starting point $\widehat{W}^{k}$ at each Newton optimization, we get the bound 
\begin{eqnarray}C \sum_{k}\big(\sup_{X}(\widehat{\mathcal{L}}_{k}(X,\mathcal{H}))-\widehat{\mathcal{L}}^{k}(\widehat{W}^{k},\mathcal{H})\big)+K \log_{2}\log_{2}(1/\epsilon)\end{eqnarray} By construction of MEMIP iterates, we also have  $\widehat{\mathcal{L}}^{k}(\widehat{W}^{k},\mathcal{H})=\widehat{\mathcal{L}}^{(k-1)}(\widehat{W}^{k},\mathcal{H})=\widehat{\mathcal{L}}^{(k-1)}(\widehat{X}^{k-1},\mathcal{H})$ where the first equality holds because for any $u$, $v$, $\widehat{W}_{uv,k}^{k}=0$ and the second because for any $u$, $v$, $j \leq k-1$, $\widehat{W}_{uv,j}^{k-1}=\widehat{X}_{uv,j}^{k-1}$. But for any $k \geq 2$, $\widehat{\mathcal{L}}^{k-1}(\widehat{X}^{k-1},\mathcal{H}) \geq \sup_{X}(\widehat{\mathcal{L}}_{k-1}(X,\mathcal{H})) - \epsilon$. Therefore the bound reformulates as \begin{eqnarray} C \sum_{k=1}^{K}\big(\sup_{X}(\widehat{\mathcal{L}}_{k}(X,\mathcal{H}))-\sup_{X}(\widehat{\mathcal{L}}_{k-1}(X,\mathcal{H})))\big) +K  (\log_{2}\log_{2}(1/\epsilon)+C\epsilon) \end{eqnarray}
which proves Proposition 4, using the notation $\sup_{X}(\widehat{\mathcal{L}}_{0}(X,\mathcal{H}))=\widehat{\mathcal{L}}_{1}(\widehat{W}^{1},\mathcal{H})$.\qed

\medskip

\noindent {\bf Remark.} The previous proposition emphasizes the key role played by the starting point $\widehat{W}^{1}$ in the speed of convergence of Newton-like methods. In our case, a good choice is for instance to select it by classical non-negative maximization techniques for objectives of type (\ref{eqn:LLhat}) (see \emph{e.g} \cite{seung2001algorithms}). Because these methods are quite fast, they can also be used for steps $k \in [2...K]$ in order to provide an alternative starting point $\widehat{W}^{k}_{+}$. The update $\widehat{X}^{k}$ is then given by either $\text{NewtonArgMax}(\widehat{\mathcal{L}}^k, \widehat{W}^{k})$ or $\text{NewtonArgMax}(\widehat{\mathcal{L}}^k, \widehat{W}^{k}_{+})$ depending on the most succesful maximization.

\section{Experimental Results}

We first evaluate MEMIP on realistic synthetic data sets. We compare it to MMEL (\cite{zhou2013learning}) and fixed exponential kernels and show that MEMIP performs significantly better in terms of prediction and triggering kernels recovery.

\subsection{Synthetic Data Sets: Experiment Setup and Results}

\noindent {\bf Data Generation} We simulate multivariate Hawkes processes by \textit{Ogata modified thinning algorithm} (see e.g. \cite {liniger2009multivariate}). Since each occurence can potentially increase stochastic rates of all events, special attention has to be paid to avoid \emph{explosion}, \emph{i.e} the occurence of an infinite number of events on a finite time window. In order to avoid such behavior, our simulated data sets verify the sufficient non-explosion condition $\rho(\Gamma)<1$ where $\rho(\Gamma)$ denotes the spectral radius of the matrix $\Gamma=(\int_{0}^{\infty}{|g_{uv}(t)}dt|)_{uv}$ (see \emph{e.g} \cite {daley2007introduction}). We perform experiments on three different simulated data sets where triggering kernels are taken as 
\begin{equation}
g_{uv}(t)=\nu_{uv}\frac{\sin \left(\frac{2\pi t}{\omega_{uv}}+\frac{\pi}{2}((u+v) \bmod 2)\right)+2}{3(t+1)^{2}}
\label{eqn:TrigKerData}
\end{equation}
We sample the periods $\omega_{uv}$ from an uniform distribution over $[1,10]$. Absolute values of normalization factors $\nu_{uv}$ are sampled uniformally from $[0,1/d[$ and their sign is sampled from a Bernoulli law of parameter $p$. Except for the toy data set, background rates $\mu_{v}$ are taken constant and sampled in $[0,0.001]$. An important feature of this choice of triggering kernels and parameters is that resulting Hawkes processes respect the aforementioned sufficient non-explosion condition. For quantitative evaluation, we simulate two quite large data sets (1) $d=300$,$p=1$ (2) $d=300$,$p=0.9$. Thus, data set (1) contains realizations of purely mutually-exciting processes whereas data set (2) has $10\%$ of inhibitive kernels. For each data set, we sample 10 sets of parameters $(\omega_{uv},\nu_{uv})_{u\geq1,v\geq1}$,$(\mu_{v})_{v\geq1}$ and simulate 400,000 i.i.d realizations of the resulting Hawkes process over $[0,20]$. The first 200,000 are taken as training set and the remaining 200,000 as test set.

\medskip

\noindent {\bf Evaluation Metrics} We evaluate the different algorithms by two metrics: (a) \emph{Diff} a normalized $L^{2}$ distance between the true and estimated triggering kernels, defined by  
\begin{eqnarray}
\mbox{Diff} = \frac{1}{d^{2}}\sum_{u=1}^{d}\sum_{v=1}^{d}\frac{\int{(\widehat{g}_{uv}-g_{uv})^{2}}}{\int{\widehat{g}_{uv}^{2}}+\int{g_{uv}^{2}}}
\end{eqnarray}
, (b) \emph{Pred} a prediction score on the test data set defined as follows. For each dimension $u \in [1...d]$ and occurence $i$ in the test set, probability for that occurence to be of type $u$ is given by $P^{true}_{i}(u)=\frac{\lambda_{u}(t_ {i})}{\sum_{v=1}^{d} \lambda_{v}(t_{i})}$. Thus, defining $AUC(d,P)$ the area under ROC curve for binary task of predicting $(1_{\left\{u_{i}=u\right\}})_{i}$ with scores $(P^{true}_{i}(d))_{i}$ and $ (P_{i}^{model}(d))_{i}$ the probabilities estimated by the evaluated model, we set

\begin{eqnarray}
\mbox{Pred}=\frac{\sum_{u=1}^{d}{(AUC(d,P^{model})-0.5)}}{\sum_{u=1}^{d}{(AUC(d,P^{true})-0.5)}}
\end{eqnarray}

\medskip

\noindent {\bf Baselines} We compare MEMIP to (a) \textbf{MMEL} for which we try various sets of number of base kernels, total number of iterations and smoothing hyperparameter, (b) \textbf{Exp} the widely used setting where $g_{uv}(t)=\nu_{uv} \exp(-\alpha t)$ and only $\nu_{uv}$ are estimated from the data. In order to give this baseline more flexibility and prediction power, we allow negative values of $\nu_{uv}$. We train three different versions with $\alpha \in \{0.1,1.0,10.0\}$.  

\medskip

\begin{figure}
\begin{center}
\includegraphics[width=1.0\linewidth]{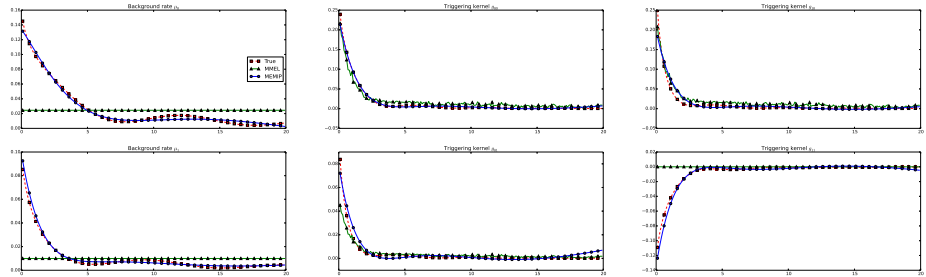}
\caption{Triggering kernels and background rates for toy data set estimated by MEMIP and MMEL algorithms vs true triggering kernels and background rate}
\label{fig:MEMIPToy3}
\end{center}
\end{figure}

\noindent {\bf Results Part 1: Visualization on a Toy Dataset} In order to demonstrate the ability of MEMIP to discover the underlying dynamics of Hawkes processes even in presence of inhibition and varying background rates, we construct the following toy bidimensional data set. Amongst the four triggering kernels, $g_{11}$ is taken negative and background rates are defined by $\mu_{0}=\frac{cos(\frac{2 \pi t}{\omega_{0}})+2}{1+t}$ and $\mu_{1}=\frac{sin(\frac{2 \pi t}{\omega_{1}})+2}{1+t}$ with parameters $\omega_{0}$ and $\omega_{1}$ sampled in $[5,15]$. We sample a set of parameters $(\omega_{uv},\nu_{uv})_{u\geq1,v\geq1}$,$(\mu_{v})_{v\geq1}$ and simulate 200,000 i.i.d realizations of the resulting Hawkes process. From Fig. \textbf{\ref{fig:MEMIPToy3}}, we observe that both compared methods MEMIP and MMEL accurately recover nonnegative triggering kernels $g_{00}$, $g_{01}$ and $g_{10}$. However, MEMIP is also able to estimate the inhibitive $g_{11}$ whereas MMEL predicts $g_{11}=0$. Varying background rates $\mu_{0}$ and $\mu_{1}$ are also well estimated by MEMIP, whereas by construction MMEL and Exp only return constant values $\bar \mu_{0}$ and $\bar \mu_{1}$. 

\medskip

\begin{figure}
\begin{center}

\includegraphics[width=1.0\linewidth]{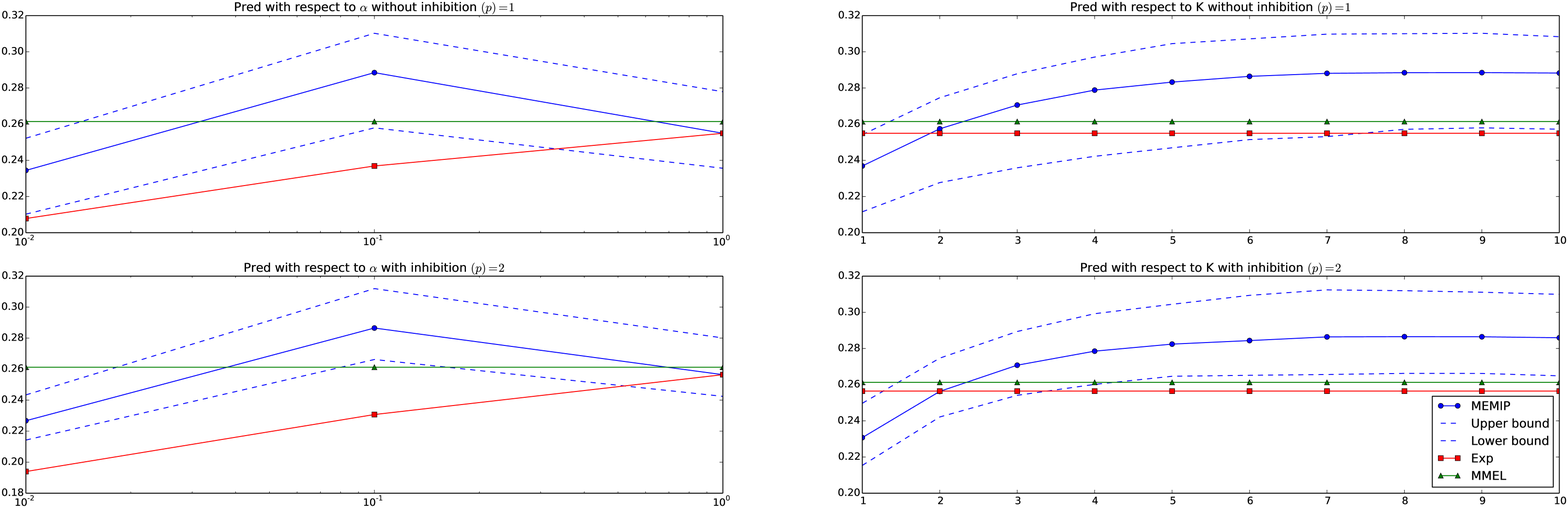}
\caption{Sensitivity to hyperparameters $\alpha$ (left) and $K$(right) for \textbf{Pred} score of MEMIP algorithm, compared to Exp and MMEL baselines on non-inhibitive simulated data set (above) and simulated data set with 10 \% inhibitive kernels (below)}
\label{fig:MEMIPSensiSimul}
\end{center}
\end{figure}

\noindent {\bf Results Part 2: Prediction Score} In order to evaluate \textbf{Pred} score of the competing methods on the generated data sets, we remove for each model the best and worst perfomance over the ten simulated processes, and average \textbf{Pred} over the eight remaining one. Empirical 10\% confidence intervals are also indicated to assess significativity of the experimental results. From Table \textbf{1}, we observe that MEMIP significantly outperforms the competing baselines for both data sets. Prediction rates are quite low for all competing methods which indicates a rather difficult prediction problem, as $90,000$ nonparametric functions are indeed to be estimated from the data. In Fig. \textbf{\ref{fig:MEMIPSensiSimul}} , we study the sensitivity of \textbf{Pred} score to $\alpha$ and $K$ for simulated data sets (1)(above) and (2)(below). Left plots show MEMIP and Exp \textbf{Pred} score with respect to $\alpha$, as well as best MMEL average score across a broad range of hyperparameters. Empirical 10\% confidence intervals are also plotted in dashed line. We see that MEMIP gives good results in a wide range of values of $\alpha$, and outperforms the exponential baseline for all values of $\alpha$. Right plots show MEMIP \textbf{Pred} score with respect to $K$ for $\alpha=0.1$, as well as best Exp and MMEL average score. We see that MEMIP achieves good prediction results for low values of $K$, and that taking $K>10$ is not necessary. For very large values of $\alpha$, we also note that MEMIP and Exp baseline are the same, because the optimal choice of $K$ for MEMIP is $K=1$.

\begin{table}
\centering
\label{table:PredSimul}
\caption{Pred score for prediction of the type of next event on simulated data sets}
\begin{tabular}{llll}
\hline
Dataset & MEMIP & MMEL & Exp \\\hline
(1) d=300,p=1 &$0.288 \in [0.258,0.310]$ &$0.261\in[0.250,0.281]$ &$0.255 \in[0.236;0.278]$ \\
(2) d=300,p=0.9 &$0.287 \in [0.266,0.312]$ &$0.261\in[0.241,0.280]$ &$0.256 \in [0.242,0.280]$ \\
\hline
\end{tabular}
\end{table}

\medskip

\begin{figure}
\begin{center}

\includegraphics[width=1.0\linewidth]{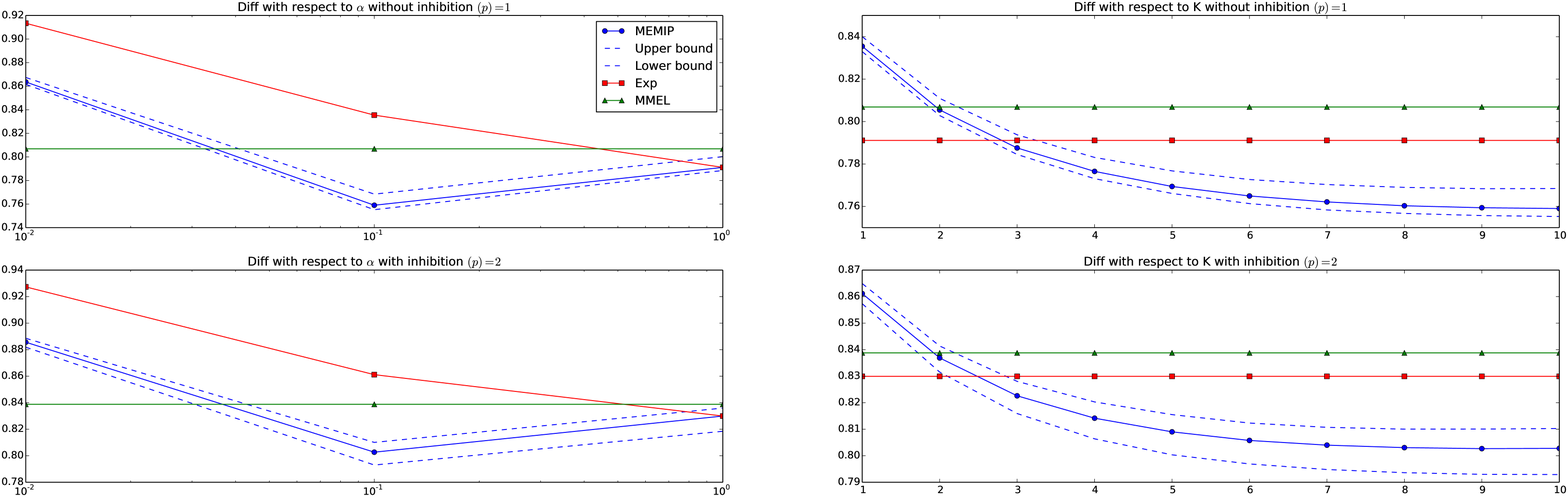}
\caption{Sensitivity to hyperparameters $\alpha$ (left) and $K$(right) for \textbf{Diff} score of MEMIP algorithm, compared to Exp and MMEL baselines on non-inhibitive simulated data set (above) and simulated data set with 10 \% inhibitive kernels (below)}
\label{fig:MEMIPSensiSimulDiff}
\end{center}
\end{figure}

\noindent {\bf Results Part 3: Accuracy of Kernel Estimation} Besides having a greater prediction power, we observe in Table \textbf{2} that MEMIP is also able to estimate the true values of triggering kernels more accurately on both data sets. In Fig. \textbf{\ref{fig:MEMIPSensiSimulDiff}}, we study the sensitivity of \textbf{Diff} score to $\alpha$ and $K$ for simulated data sets (1)(above) and (2)(below). We see that the variance of \textbf{Diff} score is very low for MEMIP, and its fitting error is significatively lower than those of other baselines at level 10\%.

\begin{table}
\centering
\label{table:DiffSimul}
\caption{Diff score for triggering kernels recovery on simulated data sets}
\begin{tabular}{llll}
\hline
Dataset & MEMIP & MMEL & Exp \\\hline
(1) d=300,p=1 &$0.759 \in [0.755,0.768]$ &$0.807 \in [0.803,0.814]$ & $0.791 \in [0.788,0.800]$ \\
(2) d=300,p=0.9 &$0.803 \in [0.793,0.810]$ &$0.839 \in [0.833,0.844]$ & $0.830 \in [0.818,0.836]$ \\
\hline
\end{tabular}
\end{table}

\medskip

\noindent {\bf Discussion} The fact that the proposed algorithm MEMIP outperforms MMEL on a non-inhibitive data set may seem surprising. Actually, even for purely mutually-exciting settings, these two algorithms can exhibit different behaviors. MMEL decomposes the triggering kernels on a low-rank set of basis functions, whereas we fix our basis functions as exponentials, in order to enjoy fast global convergence and ability to learn negative projection coefficients $X_{uv,k}$. Smoothing strategy also plays a key role in experimental results. Indeed, because the log-likelihood (\ref{eqn:linear Hawkes}) can be made arbitrarily high by the sequence of functions $(g_{uv}^{n})_{n \in N}$ defined by $g_{uv}^{n}(t) = n 1_{\left\{t \in T_{uv} \right\}}$ where $T_{uv}=\{t_{v}-t_{u} \mid (t_{u}<t_{v} \land (\exists h \in \mathcal{H} \mid (t_{v},v) \in h \land (t_{u},u) \in h ))\}$, smoothing is mandatory when learning triggering kernels by means of log-likelihood maximization. Using a $L^{2}$ roughness norm penalization $\alpha \int_{0}^{T} g'^{2}$, MMEL can face difficult dilemmas when fitting power-laws fastly decaying around $0$ : either under-estimating the rate when it is at its peak or lowering the smoothness parameter and being vulnerable to overfitting. On the contrary, MEMIP would face difficulties to perfectly fit periodic functions with a very small period, as the derivative of its order $K$ estimates can only vanish $K-1$ times.

\subsection{Experiment on the MemeTracker Data Set}
\begin{figure}
\begin{center}
\includegraphics[width=1.0\linewidth]{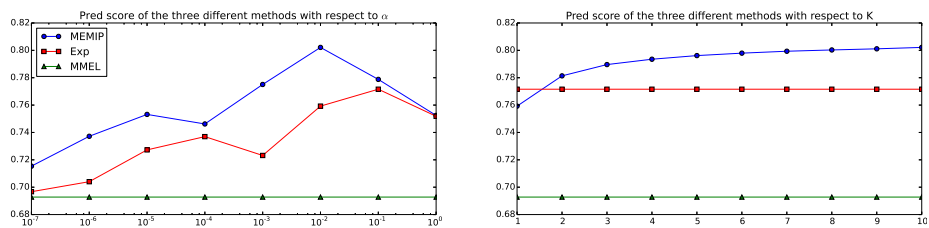}
\caption{Sensitivity to hyperparameters $\alpha$ (left) and $K$(right) for prediction score of MEMIP algorithm, compared to Exp and MMEL baselines on MemeTracker data set}
\label{fig:MEMIPSensiTracker}
\end{center}
\end{figure}

In order to show that the ability to estimate inhibitive triggering kenels and varying background rates yields better accuracy on real-world data sets, we compare the proposed method MEMIP to different baselines on the MemeTracker data set, following the experience plan exposed in \cite{zhou2013learning}. MemeTracker contains links creation between some of the most popular websites between August 2008 and April 2009. We extract link creations between the top 100 popular websites and define the occurence of an event for the $i^{th}$ website as a link creation on this website to one the 99 other websites. We then use half of the data set as training data and the other half at test data on which each baseline is evaluated by average area under ROC curve for predicting future events. From Fig. \textbf{\ref{fig:MEMIPSensiTracker}}, we observe that the proposed method MEMIP achieves a better prediction score than both baselines. Left plot shows MEMIP and Exp prediction score with respect to $\alpha$, as well as best MMEL score across a broad range of hyperparameters. We see that MEMIP gives good results in a very broad range of values of $\alpha$, and significantly outperforms the exponential baseline for all values of $\alpha$. Right plot shows MEMIP prediction score with respect to $K$ for $\alpha=0.01$, as well as best Exp and MMEL score. For $K=10$, MEMIP achieves a prediction score of $0.8021$ whereas best MMEL and Exp score are respectively $0.6928$ and $0.7716$. We note that, even for K as low as $3$, MEMIP performs the prediction task quite accurately.

\section{Conclusions}

In this paper, we propose MEMIP, which is to our knowledge the first method to learn nonparametrically triggering kernels of multivariate Hawkes processes in presence of inhibition and varying background rates. By relying on results of approximation theory, the triggering kernels are decomposed on a basis on memoryless exponential kernels. This maximization of the log-likelihood is then shown to reformulate as a concave maximization problem, that can be solved in linear complexity thanks to the Markov property verified by the proposed estimates. Experimental results on both synthetic and real-world data sets show that the proposed model is able to learn more accurately the underlying dynamics of Hawkes processes and therefore has a greater prediction power.


\bibliographystyle{splncs} 
\bibliography{LearningKernels} 

\end{document}